\magnification=1200
\baselineskip=14pt

\def\qed{{$\vrule height4pt depth0pt width4pt$}}
\def\ss{{\smallskip}}
\def\ms{{\medskip}}
\def\bs{{\bigskip}}
\def\ni{{\noindent}}

\def\s{{\sigma}}
\def\a{{\alpha}}

\def\tM{{\tilde {\cal M}}}
\def\l{{\lambda}}
\def\A{{\cal{A}}}
\def\g{{\gamma}}
\def\d{{\delta}}
\def\k{{\kappa}}
\def\b{{\beta}}
\def\Z{{\bf Z}}
\def\Q{{\bf Q}}
\def\<{{\langle}}
\def\>{{\rangle}}

\def\D{{\Delta}}
\def\T{{\bf T}}
\def\R{{\bf R}}
\def\<{{\langle}}
\def\>{{\rangle}}

\def\m1{{\quad ({\rm mod\ 1})}}
\def\L{{\Lambda}}
\def\m{{\bf m}}

\def\R1{{\cal R}_1}

\def\M{{\cal M}}
\input BoxedEPS.tex
\SetTexturesEPSFSpecial
\HideDisplacementBoxes

\centerline{\bf TORSION NUMBERS OF AUGMENTED GROUPS}
\centerline{\bf with applications to knots and links}
\smallskip

\centerline{Daniel S. Silver and Susan G. Williams} \ss
\centerline{ Revised: July, 2002}
\ms
\footnote{}{First author partially supported by \'EGIDE at CMI,
Universit\'e de Provence. Second author partially supported
by CNRS at Institut de Math\'ematiques de Luminy. Both authors 
partially supported by NSF grant DMS-0071004.}
\footnote{}{2000 {\it Mathematics Subject Classification.}  
Primary 57Q45; secondary 37B10, 11R06.}

\baselineskip=12pt
\noindent {\narrower\narrower\smallskip\noindent 
{{\bf ABSTRACT.} Torsion and Betti numbers for knots are
special cases of more general invariants $b_r$ and $\b_r$,
respectively, associated to a finitely generated group $G$ and
epimorphism $\chi:G\to \Z$. The sequence of Betti numbers is
always periodic; under mild hypotheses about
$(G,\chi)$, the sequence
$b_r$ satisfies a linear homogeneous recurrence
relation with constant coefficients. Generally, $b_r$
exhibits exponential growth rate. However, again under mild
hypotheses, the $p$-part of $b_r$ has trivial growth for any
prime $p$. Applications to branched
cover homology for knots and links are presented.}

\smallskip}
\ms
\baselineskip=14pt

\ni {\bf 1. Introduction.} A {\it knot} is a simple closed
curve in the $3$-sphere $S^3$. Knots are
{\it equivalent} if there is an orientation-preserving 
homeomorphism of $S^3$ that carries one into the other.
Equivalent knots are regarded as the same. An {\it invariant}
is a well-defined quantity that depends only on a knot
equivalence class. Two knots for which some invariant
differs are necessarily distinct. 

Associated to any knot $k$ and natural number $r$ there is a
compact, oriented $3$-manifold $M_r$, the $r$-fold cyclic
cover of $S^3$ branched over $k$. A precise definition 
can be found in [{\bf Li97}] or [{\bf Ro76}], for example.
Topological invariants of
$M_r$ are invariants of $k$. Two such invariants, the
first Betti number $\b_r$ and the order $b_r$ of the torsion
subgroup of
$H_1(M_r;\Z)$, were first considered by J.
Alexander and G. Briggs [{\bf Al28}], [{\bf AB27}]
and by O. Zariski [{\bf Za32}]. The
continuing interest in these invariants is witnessed
by numerous papers (e.g., [{\bf Go72}], [{\bf Me80}],
[{\bf We80}], [{\bf Ri90}] and [{\bf GS91}] ). We call
$b_r$ the {\it rth torsion number} of $k$. We say
that $b_r$ is {\it pure} if the corresponding Betti
number $\b_r$ vanishes (equivalently, $H_1(M_r;\Z)$
is a pure torsion group). 

Betti numbers are known to be periodic in $r$,
and they are relatively easy to compute (see Proposition
2.2). A useful formula for pure torsion numbers was
given by R. Fox in [{\bf Fo54}]. Although the proof
given by Fox was insufficient, a complete argument
was given by C. Weber [{\bf We80}]. Weber observed
that the problem of computing non-pure torsion
numbers is {\it ``$\ldots$une question plus
difficile.''}

Torsion and Betti numbers for knots are a special case of a
more general, algebraic construction that depends only on an 
{\it augmented group}, consisting of a finitely generated group
$G$ and a surjection $\chi:G \to \Z$. We define torsion and
Betti numbers in this general context. For a large class of 
augmented groups, including those that correspond to knots, we
provide a formula for all torsion numbers, generalizing the
formula of Fox. We prove that the sequence of torsion numbers
satisfies a linear recurrence relation. 


Torsion numbers tend to grow quickly as their index
$r$ becomes large. F.~Gonz\'alez-Acu\~na and H.~Short [{\bf GS91}]
and independently R.~Riley [{\bf Ri90}] proved that the sequence of
pure torsion  numbers of any knot $k$ has exponential growth
rate equal to the Mahler measure of the Alexander polynomial of $k$.
We improved upon this in [{\bf SW00}] by showing that
the entire sequence $b_r$ grows at this rate and generalizing
the result in a natural way for links. The proofs in [{\bf
SW00}] use a deep result about algebraic dynamical systems due
to D.~Lind, K.~Schmidt and T.~Ward (Theorem 21.1 of [{\bf
Sc95}]). Here we extend such results for torsion numbers
$b_r$ associated to many augmented groups. 
In contrast, we prove under suitable hypotheses that for
any prime number $p$ the {\it
$p$-component} of $b_r$ (i.e., the largest power of $p$ that divides
$b_r$) grows subexponentially. The proof relies on a
$p$-adic version of Jensen's formula, proven by G.R.~Everest and
B.N\'i Fhlath\'uin [{\bf EF96}], [{\bf Ev99}]. As a corollary we
strengthen a theorem of C. Gordon [{\bf Go72}] by proving that
for any knot the sequence of torsion numbers either is periodic
or else displays infinitely many prime numbers in the
factorization of its terms. 

In the final section we apply our techniques to the problem of
computing homology groups of branched cyclic covering spaces
associated to knots and links. 

We are grateful to Dan Flath, Adam Sikora, 
Doug Lind and Hamish Short
for useful discussions. 
The University of Maryland,
the Centre de Math\'ematituqes et Informatique 
in Marseille, and Institut de
Math\'ematiques de Luminy provided kind 
hospitality during the period
of this work. Finally, we thank the referees
for helpful comments and suggestions.
\bs


\ni {\bf 2. Augmented groups and torsion numbers.} Torsion numbers
for knots and links arise as a special case of a  general
group-theoretical quantity described below. We see that many
knot-theoretic results remain valid in the broader context.

Let $G$ be a finitely
generated group and $\chi: G \to \Z$ any epimorphism. The
pair $(G, \chi)$ is called an {\it augmented group}. Two augmented
groups, $(G_1,\chi_1)$ and $(G_2,\chi_2)$, are {\it equivalent} if
there exists an isomorphism $\phi: G_1 \to G_2$ such that
$\chi_2\circ \phi =
\chi_1$. 

For any augmented group $(G,\chi)$, the abelianization of
${\rm ker}\ \chi$ is a module $\M$ over the ring $\R1= \Z[t,
t^{-1}]$ of Laurent polynomials. Since $\R1$ is Noetherian,
$\M$ is finitely generated, expressible as
$$\M\cong\R1^N/\A\R1^M, \eqno(2.1)$$
where $\A$ is an $N\times M$-matrix over $\R1$,
for some positive integers $M, N$. By adjoining
zero columns if needed, we can assume that $M\ge
N$. 

For any
natural number
$r$, we
define $\M_r$ to be the quotient module
$$\M_r = \M/(t^r-1)\M.$$
It is clear that $\M_r$ is finitely generated as an abelian
group. Hence it decomposes as 
$$\M_r \cong \Z^{\b_r}\oplus T\M_r,$$
where $T\M_r$ denotes the torsion subgroup of $\M_r$. We define the 
{\it rth torsion number} of $(G,\chi)$ to be the order $b_r$
of
$T\M_r$. We say that $b_r$ is {\it pure} if the
Betti number
$\b_r$ vanishes. It is a straightforward matter to check
that
$\b_r$ and
$b_r$ depend only on the equivalence class of $(G, \chi)$.

The elementary
ideals $E_i$ of ${\cal M}$ form a sequence of
invariants of $(G, \chi)$. The ideal
$E_i$ is generated by the $(N-i)\times (N-i)$  minors of the matrix
$\A$ of (2.1). Since $\R1$ is a unique factorization domain, each
$E_i$ is contained in a unique minimal principal ideal; a
generator is
the $i$th {\it characteristic polynomial} $\D_i(t)$ of 
$(G,\chi)$, well defined up to multiplication by units in $\R1$. We are primarily
interested in $\D_0(t)$, which we abbreviate by $\D$.

An important class of augmented groups arises in knot
theory. For any knot $k$ in the $3$-sphere
$S^3$ the fundamental group
$G=
\pi_1(S^3-k)$ is finitely presented and has
infinite cyclic abelianization. Abelianization
provides a surjection $\chi: G \to \Z$. (More
precisely, there are two choices. The ambiguity,
which is harmless, can be eliminated by
orienting the knot.) The module
$\M$ is isomorphic to the first homology group of the infinite cyclic
cover of $S^3-k$, and it has a presentation
marix $\A$ that is square (that is, $M=N$). The
quotient module
$\M_r$ is isomorphic to the homology group $H_1(M_r, \Z)$ of the
$r$-fold cyclic cover $M_r$ of
$S^3$ branched over $k$. The $0$th characteristic polynomial $\D$
is commonly called the Alexander polynomial of $k$.  (See [{\bf
Li97}] or [{\bf Ro76}].) \bs

\ni {\bf Definition 2.1.} The {\it cyclotomic order} $\g=\g(\D)$ is
the least common multiple of those positive integers $d$ such that 
the $d$th cyclotomic polynomial $\Phi_d$ divides $\D$. If no 
cyclotomic polynomial divides $\D$ then $\g=1$. \bs

\ni {\bf Proposition 2.2.} (Cf. Theorem 4.2 of [{\bf Go72}]) For any
augmented group $(G, \chi)$ the sequence $\{\b_r\}$ of Betti
numbers satisfies  $\b_{r+\g} = \b_r$, where $\g$ is the  cyclotomic
order of
$\D$.  \bs

\ni {\bf Proof.} We adapt an argument of D.~W. Sumners
that appears in [{\bf Go72}].

Since  $\Pi ={\bf C}[t, t^{-1}]$ is a principal ideal domain, the
tensor product $\M\otimes_\Z{\bf C}$ decomposes as a direct sum
$\oplus_{i=1}^n \Pi/(\pi_i)$, for some elements $\pi_i \in
\Pi$ such that  $\pi_i|\pi_{i+1},\ 1\le i< n$. (For $0\le i< n$, the
product
$\pi_1\cdots\pi_{n-i}$ is the same as $\D_i$ up to
multiplication by units in $\Pi$.) Likewise,
$$\M_r\otimes_\Z{\bf C} \cong \oplus_{i=1}^n \Pi/(\pi_i,
t^r-1).$$
Each factor $\Pi/(\pi_i)$ can be expressed as 
$\oplus_j \Pi/((t-\a_j)^{e(\a_j)})$, where  $e(\a_j)$ are 
positive
integers, $\a_j$ ranging over the  distinct roots of $\pi_i$. Since
$$((t-\a)^{e(\a)}, t^r-1) = \cases{(t-\a) & if
$\a^r=1,$ \cr \quad \Pi &  otherwise,}$$
we see that 
$$\b_r = {\rm dim}_{\bf C} \M_r\otimes_\Z{\bf C} =
\sum_{i=1}^n l_i,$$
where $l_i$ is the number of distinct roots of $\pi_i$ that are
also $r$th roots of unity. Hence $\b_r =
\b_{(\g,r)}$, and so $\b_{r+\g} = \b_{(\g, r+\g)} = \b_{(\g,r)} =
\b_r$. \qed
\bs

In view of Proposition 2.2 it is  natural to consider a
subsequence of torsion numbers $b_{r_k}$ such that 
$\b_{r_k}$ is constant. We prove that $\{b_{r_k}\}$ is a  {\it
division sequence} in the sense that $b_{r_k}$ divides
$b_{r_l}$ whenever $r_k$ divides $r_l$.  \bs

\ni {\bf Lemma 2.3.} Assume that $\phi:{\cal N} \to {\cal N}'$
is an epimorphism of finitely generated modules over a PID. If
${\cal N}$ and ${\cal N}'$ have the same rank, then $\phi$
restricts to an epimorphism 
$\phi: T{\cal N} \to T{\cal N}'$ of torsion
submodules. \bs

\ni {\bf Proof.} It is clear that $\phi$ induces an epimorphism
$\bar \phi: {\cal N}/T{\cal N} \to {\cal N}'/T{\cal N}'$. Since
${\cal N}$ and ${\cal N}'$ have the same rank, $\bar \phi$ is an
isomorphism. If
$y \in T{\cal N}'$, then there exists an
element
$x\in {\cal N}$ such that
$\phi(x)=y$. If $x \notin T{\cal N}$, then $x$ represents a
nontrivial element of the kernel of $\bar \phi$, a
contradiction. Thus $\phi$ restricts to an epimorphism of
torsion submodules. \qed \bs

\ni {\bf Proposition 2.4.} Let $(G,\chi)$ be an augmented group.
If $b_{r_k}$ is a subsequence of torsion
numbers for which the corresponding Betti numbers $\b_{r_k}$ are
constant, then  $\{b_{r_k}\}$ is a division sequence.
\bs

\ni {\bf Proof.} If $r$ divides $s$, then clearly there
exists a surjection $\phi:\M_s\to \M_r.$ Since $\b_r =\b_s$, 
Lemma 2.3 implies that $\phi$ induces a surjection of torsion
submodules, and consequently $b_r$ divides $b_s$. \qed\bs

Given an augmented group $(G, \chi)$ such that $\M$ has a square
matrix presentation (2.1), the pure torsion numbers
$b_r$ can be computed by the following formula familiar
to knot theorists. \bs

\ni {\bf Proposition 2.5.} Assume that $(G, \chi)$ is an augmented
group such that
$\M$ has a square matrix presentation.  If $b_r$ is a pure torsion
number, then it is equal to the absolute value of 
$$\prod_{\zeta^r=1} \D(\zeta).\eqno(2.2) $$
\bs

The quantity (2.2) is equal to the resultant ${\rm Res}(\D,
t^r-1)$. In general, if $f(t)= a_0t^n+\cdots +a_{n-1}t+a_n$ and
$g(t)= b_0t^m+\cdots+b_{m-1}t+b_m$ are polynomials with integer
coefficients and zeros $\a_1,\ldots, \a_n$ and $\b_1,\ldots,
\b_m$, respectively, then the {\it resultant} of $f$ and $g$ is
$${\rm Res}(f,g)=(a_0^mb_0^n)\prod_{i,j}(\a_i-\b_j)=a_0^m\prod_i
g(\a_i)= (-1)^{mn}b_0^n\prod_j f(\b_j).$$ 
Clearly, ${\rm Res}(f_1f_2,g)={\rm Res}(f_1,g){\rm Res}(f_2,g)$ and
${\rm Res}(f,g)=(-1)^{mn}{\rm Res}(g,f)$. The resultant has an alternative
definition as the determinant of a certain matrix formed from the coefficients of
$f$ and $g$ (cf.\ [{\bf La65}]).  In particular, the resultant of integer polynomials
is always an integer.

In the case that $G$ is a knot group, formula (2.2) was 
given by R. Fox [{\bf Fo54}]. A complete proof is contained in [{\bf
We80}]. The proof of Proposition 2.5 can be fashioned along similar
lines. We will prove a more general result in Section 3. \bs

In [{\bf Le33}] D.~H. Lehmer investigated resultants
${\rm Res}(f, t^r-1)$, where $f(t)\in \Z[t]$. As he observed, it
follows from a theorem of Lagrange that the sequence $\{{\rm Res}(f,
t^r-1)\}$ satisfies a linear homogeneous recurrence relation in
$r$ with constant coefficients. 

The
general linear recurrence relation is easy to find. Assume
that
$f(t)= c_0t^d+\cdots+c_{d-1}t+c_d$ has roots
$\a_1,\ldots,\a_d$. Form the polynomials
$$\eqalign{ f_0(t)&=t-1,\cr
f_1(t)&={1 \over c_0} f(t)=\prod_{i=1}^d(t-\a_i),\cr
f_2(t)&=\prod_{i>j=1}^{d-1}(t-\a_i\a_j),\cr
&\ \ \vdots\cr
f_d(t)&=t-\a_1\a_2\cdots\a_d = t-(-1)^d{c_d \over c_0}.}$$ It is not necessary to
find the roots of $f$ in order to determine $f_0,\ldots,f_d$. The coefficients of
these polynomials are integers obtained rationally in terms of
the coefficients of $f$. Details can be found in [{\bf Le33}]. 
If $t^m+A_1t^{m-1}+\cdots+A_m$ is the
least common multiple of
$f_0,\ldots, f_d$, then ${\rm Res}(f, t^r-1)$, which we abbreviate by
$R(f,r)$, satisfies the homogeneous linear recurrence
with  characteristic polynomial $p(t)=
c_0^mt^m+c_0^{m-1}A_1t^{m-1}+\cdots+A_m$; that
is, $$c_0^m R(f,r+m)+ c_0^{m-1}A_1
R(f,r+m-1)+\cdots+A_mR(f,r)=0.
\eqno(2.3)$$ 
It is easy to see that the degree $m$ of the characteristic equation
(2.3) is not  greater than $2^d$. These facts were rediscovered by W.
Stevens [{\bf St00}]. Stevens proved that when $f$ is a reciprocal
polynomial (that is, $c_i=c_{d-i}$ for $i=0,
1,\ldots, d$) $m$ can be bounded from above by
$3^{d/2}$.

We remark that the sign of ${\rm Res}(f, t^r-1)$ is either
constant or alternating. For in the product 
$${\rm Res}(f, t^r-1)= c_0^m\prod_{i}(\a^r_i-1),$$
a pair of conjugate complex roots contributes a factor 
$(\a_i^r-1)(\bar \a_i^r-1) = |\a_i^r-1|^2$, while the 
real factors have constant or alternating sign. It follows
that $|{\rm Res}(f,t^t-1)|$ satisfies a linear recurrence
of the same order as ${\rm Res}(f, t^r-1)$; in the alternating
sign case, simply modify the characteristic polynomial by changing
the sign of alternate terms. 
\bigskip

\ni {\bf Example 2.6.} The Alexander polynomial of the
figure-eight knot (the knot $4_1$ in tables) is $\D(t)=t^2-3t+1$.
Since neither root has modulus one, all of the  torsion numbers
of $k$ are pure. The polynomials $f_i$ are $f_0(t)=f_2(t)=t-1$
and $f_1(t)=\D(t).$ The least common multiple is
$t^3-4t^2+4t-1$, and hence $b_r$ satisfies:
$b_{r+3}-4b_{r+2}+4b_{r+1}-b_r=0$. Using the initial conditions
$b_0=0, b_1=1, b_2=5$, other values can now be quickly
computed. 

The  torsion numbers for the figure-eight knot produce some
surprisingly large prime factors. According to calculations done
with Maple, $b_{1361}$ is the square of a prime with 285 digits.

Lehmer, who considered this example in [{\bf Le33}], albeit for
much smaller values of $r$, was interested in producing new prime
numbers. He observed that the factors of $R(f,r)$ satisfy a severe
arithmetical constraint, and he proposed that if $R(f,r)$ grows
with a relatively small exponential  growth rate, then these numbers
will likely display large prime factors. Lehmer did not give any
proof of the  assertion about prime factors, but rather used it
heuristically. A survey of Lehmer's efforts together
with new results in these directions can be found in [{\bf EEW00}].
\bs

\ni {\bf Definition 2.7.} Assume that
$$f(t)= c_0t^d+\cdots+c_{d-1}t+c_d =
c_0\prod_{i=1}^d(t-\a_i)$$ is a polynomial with
complex coefficients, $c_0\ne 0$. The {\it Mahler measure} of
$f$ is
$$M(f)= |c_0|\prod_{i=1}^d\max\{1,|\a_i|\}.$$
The empty product is assumed to be $1$, so that the Mahler measure
of a nonzero constant polynomial $f(t)=c_0$ is $|c_0|$. By
convention, the Mahler measure of the zero polynomial is zero.

Clearly, Mahler measure is multiplicative; that is,
$M(fg)=M(f)M(g)$, for $f,g\in {\bf C}[t]$.  The following is proved
in [{\bf GS91}] and [{\bf Ri90}]. We sketch the 
argument.
\bs

\ni {\bf Proposition 2.8.} Let
$f$ be a polynomial with integer coefficients. The
subsequence $R(f,r_k)$ of nonvanishing resultants  has
exponential growth rate $M(f)$; that is,
$$\lim_{r_k\to \infty\hfil}|{\rm Res}(f,t^{r_k}-1)|^{1/{r_k}} =
M(f).$$\bs

\ni {\bf Sketch of Proof.} Let $f(t)=
c_0t^d+\cdots + c_{d-1}t+c_d$. Assume that $c_0
\ne 0$ and that 
$\a_1, \ldots, \a_d$ (not necessarily distinct) are the roots
of $f$. Then 
$$|{\rm Res}(f, t^r-1)|^{1/r} =
|c_0|\prod_{i=1}^d|\a_i^r-1|^{1/r}.$$
The condition that the resultant does not vanish is equivalent
to the statement that no root $\a_i$ is an $r$th root of
unity. Consider the subsequence of natural integers $r$ for 
which this is the case. Note that if $|\a_i|<1$, then 
the factor $|\a_i^r-1|^{1/r}$ converges to 1 as $r$ goes to
infinity. On the other hand, if $|\a_i|>1$, then for
sufficiently large $r$ we have
$${1 \over 2}|\a_i|^r \le |\a_i|^r-1\le |\a_i^r-1|\le
|\a_i|^r+1 \le 2|\a_i|^r.$$ Taking $r$th roots we see that
$|\a_i^r-1|^{1/r}$ converges to $|\a_i|$.  

When some root $\a_i$ lies
on the unit circle the nonzero values of
$|\a_i^r-1|$ can fluctuate wildly. In this
case the analysis is more subtle. 
Gonz\'alez-Acu\~na and Short use results
of A. Baker [{\bf Ba77}] and A.O. Gelfond [{\bf Ge35}]
to obtain estimates. In [{\bf GS91}] it is shown that
if
$|\a_i^r|\ne 1$, then
$$C\exp \{-(\log r)^6\} < |\a_i^r-1|\le 2,$$
where $C$ is a positive constant that depends only on $f$.
As in the case that $|\a_i|<1$ we have that
$|\a_i^r-1|^{1/r}$ converges to $1$. 

The conclusion of Proposition 2.8 follows. \qed
\bs

The following is immediate from Proposition 2.8.\bs

\ni {\bf Corollary 2.9.} Assume that $(G,\chi)$ is an augmented group
for which the module $\M$ has a square matrix presentation.  Then the
subsequence of $\{b_r\}$ consisting of pure torsion numbers has
exponential growth rate equal to $M(\D)$. \bs

We can extend the
conclusion of Proposition 2.8 to the
entire sequence of resultants by using results
from the theory of
algebraic dynamical systems. Only the
essential elements of the theory are sketched
below. Readers unfamiliar with dynamical
systems might refer to [{\bf EW99}]. 

In brief, to an augmented group $(G,\chi)$ we
associate a compact space and a
homeomorphism $\s$ from the space to itself. 
The
fixed points of  $\s^r$ form a closed subspace
consisting of exactly $b_r$ connected components.
Topological techniques are available to compute the
exponential growth rate of $b_r$, and it coincides
with $M(\D)$. \bs

\ni {\bf Theorem 2.10.} Assume that $(G,\chi)$ is an augmented group
such that either (i) the module $\M$ has a square presentation
matrix; or (ii) $\M$ is torsion-free as
an abelian group. Then the sequence $\{b_r\}$ of torsion numbers has
exponential growth rate equal to $M(\D)$.  \bs

\ni {\bf Proof.}  Let $\M^\wedge$ denote the Pontryagin dual
${\rm Hom}(\M,\T)$; that is, the topological group of homomorphisms
$\rho$ from $\M$ to the additive circle group 
$\T={\bf R}/\Z$.
Here $\M$ has
the discrete topology, and $\M^\wedge$ the
compact-open topology.  
Multiplication by $t$ in $\M$ induces a homeomorhism $\s$ of 
$\M^\wedge$ defined by $\s(\rho)(a)= \rho(t a)$, for
any
$\rho\in \M^\wedge$ and all $a\in\M$. 
The dual of $\M_r=\M/(t^r-1)\M$ is the subspace ${\rm Fix}(\s^r)=
\{\rho\in \M^\wedge\mid \s^r\rho=\rho\}$, the set of points of 
$\M^\wedge$ with {\it period} $r$. 

Since $\M_r = \Z^{\b_r}\oplus T\M_r$, the dual $\M_r^\wedge$ is
homeomorphic to $\T^{\b_r}\times T\M_r.$ This follows from
two facts: $\Z^\wedge$ is isomorphic to $\T$; and 
$A^\wedge$ is isomorphic to $A$ for any finite abelian
group. Hence the number of connected
components of $\M_r^\wedge$ is equal to the cardinality of $T\M_r$,
which by definition is the torsion number $b_r$. Each
component is a torus of dimension $\b_r$, a beautiful
fact but one that we will not use here.

The number of connected components of $\M_r^\wedge$ is the
same as the number $N_r$ of connected components of
${\rm Fix}(\s^r)$. Theorem 21.1(3) of [{\bf Sc95}] states that the
exponential growth rate of
$N_r$ is equal to the topological entropy of $\s$.
(The proof of this deep result uses a definition of
topological entropy in terms of {\it separating sets}.
For an elementary 
discussion of the theorem see [{\bf EW99}].) Further,
if
$\M$ has a presentation (2.1) with square matrix $\A$,
then the topological entropy of
$\s$ is equal to $M(\D)$. 
(See Example 18.7(1) of [{\bf Sc95}]). Thus if the hypothesis
(i) is satisfied, then we are done. 

If $\M$ is torsion-free as an abelian group, then 
again
the topological entropy of $\s$ is equal to $M(\D)$
by Lemma 17.6 of [{\bf Sc95}]. \qed\bs

The hypotheses of Theorem 2.10 cannot be dropped, as the
following example illustrates. \bs

\ni {\bf Example 2.11.}  Consider the augmented group $(G,\chi)$
such that 
$$G= \langle x,a \mid x^{-2}a^2xa^{-6}xa^2,\
x^{-3}axa^{-4}xa^4xa^{-1}\rangle,$$
and $\chi:G \to \Z$ maps $x\mapsto 1$ and $a\mapsto 0.$ A
straightforward calculation shows that $\M\cong \R1/( 2f,
(t-1)f),$ where $f(t)=t^2-3t+1$. The Alexander polynomial
$\D$ is ${\rm g.c.d.}(2f,(t-1)f) = f$, and it has Mahler measure
greater than $1$. However, the topological entropy of the
homeomorphism
$\s$ is zero by Corollary 18.5 of [{\bf Sc95}]. As in the proof of 
the theorem above, it follows that the torsion numbers
$b_r$ have trivial exponential growth rate; that is, $\limsup_{r\to
\infty} b_r^{1/r}=1$. 
\bs


\ni {\bf 3. Extended Fox's formula and recurrence.} Let $(G,
\chi)$ be an augmented group, and $\A$ the $N\times M$
presentation matrix for the $\R1$-module $\M$ as in (2.1).
For any positive integer $r$ we can obtain a presentation matrix
for the finitely generated abelian group $\M_r$ by replacing
each entry $q(t)$ of $\A$ by the $r\times r$ block $q(C_r)$, where
$C_r$ is the companion matrix of $t^r-1$, 
$$C_r = \pmatrix{0&1&0&\cdots&0\cr
0&0&1&\cdots&0\cr
\vdots&&&&\vdots\cr
0&0&0&\cdots&1\cr
1&0&0&\cdots&0}.$$We
call the resulting $rN\times rM$ matrix $\A(C_r)$.  The
proof is not difficult. The torsion number $b_r$ is equal to the
absolute value of the product of the nonzero elementary divisors of
$\A(C_r)$.

Assume first that $\M$ is a cyclic module. Then 
$\A$ is the $1\times1$ matrix $(\D(t))$, and the 
$r\times r$
matrix $(\D(C_r))$ presents
$\M_r$.  The Betti number $\b_r$
is the number of zeros of $\D$ that are
$r$th roots of unity.
When it vanishes the matrix $(\D(C_r))$ is nonsingular. Then all
elementary divisors of the matrix are nonzero, and their product is
equal (up to sign) to the product of the eigenvalues, which is the
determinant.  Fox's formula (Proposition 2.5) follows
by choosing a basis for
${\bf C}^r$ that diagonalizes the companion matrix $C_r$; we then see
that the eigenvalues of $\D(C_r)$ are $\D(\zeta)$, where $\zeta$
ranges over the $r$th roots of unity. In general,
$\b_r$ is equal to 
$$s=\sum_{ d|r \atop \Phi_d \vert \D}{\rm deg}\ \Phi_d=\sum_{ d|r
\atop
\Phi_d
\vert
\D}\phi(d),$$ where $\Phi_d$ is, as before, the $d\,$th cyclotomic polynomial,
and $\phi$ is  Euler's phi function. We appeal to the following result, a special
case of Theorem 2.1 of [{\bf MM82}].
\bs

\ni {\bf Lemma 3.1.} Let $A$ be an integral $r\times r$ matrix
with rank $r-s$. Suppose that $R$ is an integral $s\times r$ matrix
with an $s\times s$ minor invertible over $\Z$ such that 
$RA=0$ and $AR^T=0$ (where $R^T$ denotes the
transpose matrix). Then the product of the nonzero
 eigenvalues of $A$ is equal to $\pm{\rm det}(RR^T)$
times the product of the nonzero elementary divisors
of $A$. \bs

\ni {\bf Example 3.2.} Suppose that we have a factorization $t^r-1=
\Phi\cdot \Psi$ in
$\Z[t]$.  Set $A=\Phi(C_r)$. Then we can construct a matrix $R$
satisfying the hypotheses of Lemma 3.1. We regard $\R1/(t^r-1)$ as a
free abelian group with generators $1,t, \ldots, t^{r-1}$. Then the
rows of 
$A$ represent the polynomials $\Phi, t \Phi, \ldots, t^{r-1} \Phi$
(modulo $t^r-1$). The rank of $A$ is $r-s$, where
$s= {\rm deg}\ \Phi$. We take 
$R$ to be the $s\times r$ matrix with rows representing $\Psi, t\Psi,
\ldots, t^{s-1} \Psi$. Consider first the product
$RA$.
Regarding the product of the $i$th row of $R$ with
$A$ as a linear combination of the rows of $A$, we
see that it represents the polynomial
$t^{i-1}\Psi\cdot \Phi\equiv 0$ (modulo $t^r-1$).
Hence $RA=0$. 

The columns of $A$ represent
the polynomials $\Phi(t^{-1}),  t\Phi(t^{-1}), 
\ldots, t^{r-1}\Phi(t^{-1})$,
and so the $i$th column of $AR^T$ represents
$\Phi(t^{-1})\cdot t^i\Psi(t)$ (modulo
$t^r-1$). Since $\Phi$ is a product of
cyclotomic polynomials, we have $t^{{\rm deg}\
\Phi}\Phi(t^{-1})= \pm \Phi(t)$. (A cyclotomic
polynomial has this property since its set of 
roots is preserved by inversion, and its 
leading and constant
coefficients are $\pm 1$.) So
$AR^T$ is also zero. 

Since
the degree of
$t^i
\Psi$ is less than $r$ for
$i<s$, the $s\times s$ minor consisting of the first
$s$ columns of $R$ is upper triangular. The diagonal
entries are the constant term of $\Psi$, which must be
$\pm 1$. Hence this minor is invertible over $\Z$.

The matrix $\A$ presents $\R1/(\Phi, t^r-1) \cong \R1/(\Phi)$, a
free abelian group, so the product of its
elementary divisors is $1$.  Lemma 3.1 implies
that
${\rm det}(RR^T)$ is equal up to sign to the product of the
nonzero eigenvalues of
$\Phi(C_r)$; that is,
$${\rm det}(RR^T)=\pm\prod_{\zeta^r=1\atop \Phi(\zeta)\ne0}
\Phi(\zeta).\eqno(3.1)$$ \bs

%
%

\ni {\bf Theorem 3.4.} Assume that $(G,\chi)$ is an augmented group
such that $\M\cong \R1/(\D)$. For any positive integer $r$, let $\Phi$
be the product of the distinct cyclotomic polynomials $\Phi_d$ such
that $d\vert r$ and $\Phi_d\vert \D$. Then 
$$b_r=\biggr\vert\prod_{\zeta^r=1\atop\D(\zeta)\ne0}\biggr({\D\over\Phi}\biggr)(\zeta)\biggr\vert.\eqno(3.2)$$

\bs
\ni {\bf Remarks 3.5.} (i) We follow the convention that if no
cyclotomic polynomial divides $\D$, then $\Phi=1$. Clearly $b_r$ is a
pure torsion number if and only if $\Phi=1$. In this case (3.2)
reduces to Fox's formula (2.2). \ss
\quad (ii) See [{\bf Sa95}] and [{\bf HS97}] for more calculations
and estimations of torsion numbers $b_r$ arising from link groups.

\bs

\ni {\bf Proof of Theorem 3.4.} We write $\D$ as $\Phi\cdot g$,  for
some $g\in \Z[t]$. The matrix $\D(C_r)$, which presents
$\M_r = \R1/(\D, t^r-1)$, has rank $r-{\rm deg}\ \Phi$.
The rank is the same as that of
$\Phi(C_r)$. Consider the matrix $R$ of Example 3.2. We have
$R\D(C_r)= (R\Phi(C_r))g(C_r)=0$ and also $\D(C_r)R^T
= (\Phi(C_r)g(C_r))R^T=g(C_r)(\Phi(C_r)R^T)=0.$
Formula (3.2) now follows from Lemma 3.1 together
with (3.1).
\qed \bs 

If $(G,\chi)$ is an augmented group such that $\M$ is a direct
sum of cyclic modules, then Theorem 3.4 can be applied to each
summand and the terms produced by (3.2) multiplied together in order
to compute $b_r$. 

When $\M$ is not necessarily a direct sum of cyclic modules, but it
is torsion-free as an abelian group, then it is ``virtually'' a direct
sum of cyclic modules by the following lemma, which appears as Lemma
9.1 in [{\bf Sc95}]. The main idea of the proof is to consider the
natural injection of $\M\hookrightarrow
\M\otimes_\Z\Q$, and use the fact that 
$\M\otimes_\Z\Q$ is a finitely generated module over the ring
$\Q[t^{\pm 1}]$, which is a principal ideal domain. 

We recall that a polynomial in $\Z[t]$  is said  to be
primitive if the only constants that divide it are
$\pm 1$.  \bs

\ni {\bf Lemma 3.6.} Assume that $\M$ is a finitely generated
$\R1$-module that is torsion-free as an abelian group. Then there
exist primitive polynomials $\pi_1, \ldots, \pi_n\in\Z[t]$ such that 
$\pi_i\vert \pi_{i+1}$ for all $i=1, \ldots, n-1$, and an $\R1$-module
injection $i:\M \to \M' = \R1/(\pi_1)\oplus\cdots\oplus
\R1/(\pi_n)$ such that $\M'/i(\M)$ is finite. \bs

For notational convenience we identify $\M$ with its image in
$\M'$.  Consider the mappings $\mu:\M\to \M$ and $\mu':\M' \to
\M'$ given by
$a\mapsto (t^r-1)a$. Clearly ${\rm ker}\ \mu$ is a submodule of
${\rm ker}\ \mu'$.  We define $\k(r)$ to be the index $\vert
{\rm ker}\ \mu':{\rm ker}\ \mu\vert$. Let $b_r'$ denote the order of the torsion
subgroup of $\M'/(t^r-1)\M'$. The proof of the following theorem
extends techniques of [{\bf We80}].
\bs

\ni {\bf Theorem 3.7.} Let $(G,\chi)$ be an augmented group.
If $\M$ is torsion-free as an abelian group, then for any positive
integer $r$, 
$$b_r = {{b'_r}\over{\k(r)}}.\eqno(3.3)$$
Moreover, if $\g$ is the cyclotomic order of $\D$, then
$\k(r+\g)=\k(r)$ for all $r$. \bs

\ni {\bf Lemma 3.8.} Let 
$$0 \to A_1 \to A_2 \to \cdots \to A_m \to 0$$
be an exact sequence of finite abelian groups. Then 
$$\prod|A_{\rm even}|= \prod|A_{\rm odd}|.$$ \bs

Lemma 3.8 is easily proved using induction on $m$. We leave the
details to the reader. \bs

\ni {\bf Proof of Theorem 3.7.} Consider the finite quotient $p:\M'
\to
\M'/\M$ and  mapping $\bar \mu: \M'/\M \to \M'/\M$ given by $a
\mapsto (t^r-1)a$. The exact diagram
$$\matrix{0&\to &\ \ \M &{\buildrel i \over \to}&\ \ \  \M'&{\buildrel
p \over \to}&\ \ \M'/\M&\to&0\cr
&&\mu\downarrow&&\mu'\downarrow&&\bar\mu \downarrow\cr
0&\to &\ \ \M &{\buildrel i \over \to}&\ \ \  \M'&{\buildrel
p \over \to}&\ \ \M'/\M&\to&0}$$
induces a second exact diagram 
$$\matrix{&&0&&0&&0&&\cr
&&\big\downarrow &&\big\downarrow&&\big\downarrow&& \cr
0&\to &{\rm ker}\ \mu &{\buildrel i \over
\to}&{\rm ker}\ \mu'&{\buildrel p \over \to}&{\rm ker}\ \bar\mu&\to&0 \cr
&&\big\downarrow
&&\big\downarrow&&\big\downarrow&& \cr
0&\to &\M &{\buildrel i \over \to}& \M'&{\buildrel
p \over \to}&\M'/\M&\to&0\cr
&&\hskip-9pt\mu\big\downarrow
&&\hskip-9pt\mu'\big\downarrow&&\hskip-9pt\bar\mu\big\downarrow&& \cr
0&\to &\M &{\buildrel i \over \to}& \M'&{\buildrel p \over
\to}&\M'/\M&\to&0\cr && \big\downarrow &&
\big\downarrow&&\big\downarrow&&
\cr 0&\to & \M_r &{\buildrel \bar i \over \to}&
\M'_r&{\buildrel \bar p \over \to}& {\rm coker}\ \ \bar\mu&\to&0\cr
&&\big\downarrow&&\big\downarrow&&\big\downarrow \cr
&&0&&0&&0&&\cr}$$
and hence by the Snake Lemma we obtain a long exact sequence
$$0\to {\rm ker}\ \mu\ {\buildrel i \over \to}\ {\rm ker}\ \mu'\ {\buildrel p\over
\to}\ {\rm ker}\ \bar \mu\ {\buildrel d \over \to}\ \M_r {\buildrel
\bar i\over \to}\ \M'_r\ {\buildrel \bar p \over \to}\ {\rm coker}\ \
\bar \mu \to 0. \eqno(3.4)$$
Let $T\M_r$ and $T\M'_r$ be the torsion subgroups of $\M_r$ and
$\M'_r$, respectively. Since ${\rm ker}\ \bar \mu$ is finite, its
image under the connecting homomorphism $d$ is contained in $T\M_r$.
Also,
$\bar i$ maps $T\M_r$ into $T\M'_r$. Hence we have an induced
sequence
$$0\to {\rm ker}\ \mu\ {\buildrel i \over \to}\ {\rm ker}\ \mu'\ {\buildrel p\over
\to}\ {\rm ker}\ \bar \mu\ {\buildrel d \over \to}\ T\M_r {\buildrel \bar
i\over \to}\ T\M'_r\ {\buildrel \bar p \over \to}\ {\rm coker}\ \ \bar \mu
\to 0. \eqno(3.5)$$
It is not difficult to verify that (3.5) is exact. The only nonobvious
thing to check is that the kernel of $\bar p$ is contained in the
image of $\bar i$. To see this, assume that $\bar p(y)=0$. By the
exactness of (3.4) there exists an element $x\in \M_r$ such that 
$\bar i (x)=y$. If $x \notin T\M_r$, then the multiples of $x$ are
distinct in $\M_r$ and each maps by $\bar i$ into the finite group
$T\M'_r$, contradicting the fact that ${\rm ker}\ \bar i = d({\rm
ker}\ \bar
\mu)$ is finite.

The following sequence is exact.
$$0 \to {\rm ker}\ \mu'/i({\rm ker}\ \mu) \to {\rm ker}\ \bar \mu\
\to\ T\M_r \to T\M'_r \to {\rm coker}\ \ \bar \mu \to 0.\eqno(3.6)$$ 

Since $\M'_r/\M_r$ is finite, ${\rm ker}\ \bar \mu$ and ${\rm coker}\
\ \bar \mu$ have the same order. Lemma 3.8 now completes the proof of
(3.3),  $\k(r)$ being the order of  ${\rm ker}\ \mu'/i({\rm ker}\ \mu)$.

The modules $\M$ and $\M'$ have characteristic 
polynomial $\pi_n$. 
Since $\M$ embeds in $\M'$ with finite index, a prime
polynomial annihilates a nonzero element of
$\M$ if and only if it annihilates a nonzero 
element of $\M'$. Such polynomials are exactly the
prime divisors of $\pi_n$.  It follows that ${\rm
ker}\ \mu$ and ${\rm ker}\ \mu'$ are
both  periodic, with period equal to the least common
multiple $\g$ of the positive integers $d$ such that
$\Phi_d$ divides $\D$. Hence the same is true for
$\k(r)$. \qed
\bs

\ni {\bf Theorem 3.9.} Assume that $(G,\chi)$ is an augmented group
such that $\M$ is a direct sum of cyclic modules or is torsion
free as an abelian group. Then the set of torsion numbers $b_r$
satisfies a linear homogeneous recurrence relation with constant
coefficients. \bs

\ni {\bf Proof.} Write 
$$\D= \biggr(\prod_{d\in D}\Phi_d^{e_d}\biggr)\cdot g,$$
where $D = \{d\ :\ \Phi_d|\D\,\}$, and let $\g$ be the cyclotomic order of $\D$. 
We will
show that for each $R\in \{0,\ldots,\g-1\}$, the subsequence of $b_r$ with
$r$ congruent to $R$ modulo $\g$ satisfies 
$$b_r=C_R\,r^{M_R}|{\rm Res}(g, t^r-1)|,\eqno(3.7)$$
where $C_R$, $M_R$ are constants, 
$$M_R=\sum_{d\in D \atop d|R}\phi(d)(e_d-1)\leq M=\sum_{d\in D }\phi(d)(e_d-1).$$
As we saw in section 2, the sequence $|{\rm Res}(g, t^r-1)|$ satisfies a linear
homogeneous recurrence relation with characteristic polynomial $p$ of degree at most
$2^{\deg g}$. We may normalize $p$ to be monic,
$p(t)=\prod_j (t-\l_j)^{n_j}$, with $\l_j$ distinct.
The general solution to this recurrence relation has
the form $\sum_jq_j(r)\l_j^r$, where $q_j$ is a
polynomial of degree less than $n_j$ (see [{\bf
Br92}], Theorem 7.2.2, for example.) Each of the
sequences
$a_r^{(R)}= C_Rr^{M_R}|{\rm Res}(g, t^r-1)|$ satisfies
the recurrence relation given by $\hat p (t)
=\prod_j(t-\l_j)^{n_j+M}$. It also satisfies  the recurrence
relation given by $P(t)= \prod_j(t^\g-\l_j^\g)^{n_j+M}$, since
$\hat p$ divides
$P$. Because the powers of $t$ occurring
in $P$ are all multiples of $\g$, the latter
recurrence relation also describes the sequence
$\{b_r\}$, which is composed of the subsequences
$b_{R+\g n}= a_{R+\g n}^{(R)}$.  
We note that the degree of $P$ is at most $\g
(M+1)2^{{\rm deg}\ g}$.

First we consider the case
when
$\M$ is cyclic. Given $R$ we set
$$\Phi=\prod_{d\in D\atop d|R}\Phi_d.$$
By Theorem 3.4 we have 
$$\eqalign{b_r&=\biggr\vert\prod_{\zeta^r=1\atop\D(\zeta)\ne0}\biggr( {\D \over
\Phi}\biggr)(\zeta)\biggr\vert =\bigg|{\rm Res}\bigg({\D\over\Phi}, {t^r-1\over
\Phi}\bigg)\bigg|\cr &=
\prod_{d\in D}\bigg|{\rm Res}\bigg(\Phi_d, {t^r-1\over
\Phi}\bigg)\biggr\vert^{e'_d}\biggr\vert{\rm Res}\bigg(g, {t^r-1\over
\Phi}\bigg)\biggr\vert,\cr}$$  where 
$$e'_d=\cases{e_d-1&if $d|R$,\cr
                 e_d&if $d\not\vert R$.}$$
For each
$d$ dividing $R$,
$$\eqalign{{\rm Res}(\Phi_d,
{t^r-1 \over \Phi}) &= \prod_{\Phi_d(\omega)=0}{t^r-1\over
\Phi(t)}\biggr\vert_{t=\omega}\cr
&=\prod_{\Phi_d(\omega)=0}{(t^d-1)(1+t^d+
\ldots+t^{(r/d-1)d})\over
\Phi_d(t)\hat\Phi(t)}\biggr\vert_{t=\omega}\quad 
({\rm where}\ \hat\Phi= \Phi/\Phi_d)\cr &=\prod_{\Phi_d(\omega)=0}\biggr[
{t^d-1\over
\Phi_d(t)}\biggr\vert_{t=\omega}\cdot
{r/d\over
\hat\Phi(\omega)}\biggr]=C_d\cdot r^{\phi(d)},}$$
where $C_d$ depends only on $d$ and $R$. 
For $d\in D$ not dividing $R$,
$${\rm Res}(\Phi_d,
{t^r-1 \over \Phi}) = \prod_{\Phi_d(\omega)=0}{\omega^r-1\over
\Phi(\omega)}$$
is constant for $r$ congruent to $R$ modulo $\g$,
since $d$ divides $\g$. Finally, 
$${\rm Res}(g,
{t^r-1\over
\Phi})= c_0^{r-{\rm deg}\
\Phi}\prod_{g(\a)=0}{\a^r-1\over
\Phi(\a)},$$
where $c_0$ is the leading coefficient of $g$;
the expression can be rewritten as 
$C\,{\rm Res}(g, t^r-1),$ where $C$ depends
only on $R$. Thus  we can express
$b_r$ in the desired form (3.7) for all $r$ congruent to $R$ modulo $\g$.

For the case when $\M$ is a direct sum of cyclic modules $\R1/(\pi_1)\oplus\cdots\oplus
\R1/(\pi_n)$ we apply the above result to each summand and use the 
facts that
$\D=\pi_1\ldots\pi_n$ and $b_r$ is the product of the torsion numbers of the
summands to see that equation (3.7) still holds. Finally, if
$\M$ is torsion free as an abelian group, we use Theorem 3.7. \qed
\bs

\ni {\bf 4. Prime parts of torsion numbers.} 
We recall
Jensen's formula, a short argument for which can be found in 
[{\bf Yo86}].
\bs

\ni {\bf Lemma 4.1.} [Jensen's formula] For any complex number $\a$,
$$\int_0^1\log|\a-e^{2\pi i \theta}|d\theta = \log\max\{1,|\a|\}.$$
\bs

By Lemma 4.1 the Mahler measure $M(f)$
of a nonzero polynomial with complex coefficients can be computed as
$$\exp \int_0^1\log|f(e^{2\pi i \theta})|d\theta.$$
This observation motivated the definition of Mahler measure for
polynomials in several variables. (See [{\bf Bo81}] or [{\bf EW99}],
for example.) 

In [{\bf EF96}], [{\bf Ev99}] G.R.~Everest and B.N\'i Fhlath\'uin
proved a $p$-adic analogue of Jensen's formula, which we describe.
Assume that $\a$ is an algebraic integer lying in a finite extension
$K$ of {\bf Q}. For every prime  $p$ there is a $p$-adic absolute
value $|\cdot |_p$, the usual Archimedean absolute value
corresponding to $\infty$. We recall the
definition  (see [{\bf La65}] for
more details): If
$p$ is a prime number, then
$|p^rm/n|_p=1/p^r$, where $r$ is an integer, and $m,n$
are nonzero integers that are not divisible by $p$. 
By convention,
$|0|_p = 0$. Each
$|\cdot |_p$ extends to an absolute value $|\cdot |_v$
on
$K$. Let
$\Omega_v$ denote the smallest field which is
algebraically closed and complete with respect to 
$|\cdot |_v$. Let $\T_v$  denote the closure of the group of all
roots of unity, which is in general locally compact. Note that if
$p=\infty$, then $\Omega_v={\bf C}$ and $\T_v=\T$.  Everest and 
Fhlath\'uin define 
$$M_{\T_v}(t-\a)=\exp \int_{\T_v}\log|t-\a|_vd\mu = \exp \lim_{r\to
\infty}{1\over r}\sum_{\zeta^r=1}\log|\zeta-\a|_v.$$
Here $\int$ denotes the Shnirelman integral, given by the limit of
sums at the right, where one skips over the undefined summands. The
above integral exists even if $\a\in
\T_v$, in which case it can be shown to be zero. Moreover, one has 
$$\int_{\T_v}\log |t-\a|_vd\mu = \log\max\{1,|\a|_v\},\eqno(4.1)$$
which Everest and Fhlath\'uin refer to as a $p$-adic analogue of 
Jensen's formula. 

Recall that $b_r^{(p)}$ denotes the  {\it $p$-component} of
$b_r$, the largest power of $p$ that divides $b_r$. The {\it
content}  of $f\in {\bf Z}[t]$ is the greatest common divisor of the coefficients.
Using (4.1) we will prove\bs

\ni {\bf Theorem 4.2.} Let $(G,\chi)$ be an augmented group, and let
$p$ be a prime. 
\ss

\item{}(i) If $\M$ has a square matrix presentation and
$\D(t)\neq 0$, then the sequence
$\{b_{r_k}\}$ of pure torsion numbers satisfies
$$\lim_{{r_k}\to \infty}({b_{r_k}^{(p)})}^{1/{r_k}} = ({\rm
content}\ \Delta)^{(p)}.$$
\ss
\item{}(ii) If $\M$ is a direct sum of cyclic modules, then the
sequence of all torsion numbers satisfies
$$\lim_{r\to \infty}({b_r^{(p)})}^{1/r} = ({\rm
content}\ \Delta)^{(p)}.$$\ss

(iii) If $\M$ is torsion free as an abelian group, then 
$$\lim_{r\to \infty}({b_r^{(p)})}^{1/r}=1.$$
\bs

\ni {\bf Example 4.3.} For any positive integer $m$, consider the
augmented group
$(G,\chi)$ such that $G=\langle x,y \mid y^mx = xy^m\rangle$ and 
$\chi:G\to \Z$ maps $x\mapsto 1$ and $y\mapsto 0$. One verifies
that $\M\cong \R1/(m(t-1))$. The quotient module
$\M_r$ is isomorphic to $\Z^r/A_r\Z^r$, where
$$A_r = \pmatrix {m &0&0 &0&\cdots&-m \cr -m&\ m&0&\cdots&&\ \
0\cr 0&-m&m&0&\cdots&\ \ 0\cr &&\vdots\cr 0&0&&\cdots
&-m&\ \ m\cr}.$$ 
The matrix is equivalent by elementary row
and column operations to the diagonal matrix
$$\pmatrix{m\cr &\ddots \cr&&&m\cr&&&&0}.$$
Hence $\M_r\cong \Z\oplus(\Z/m)^{r-1}$, and so $b_r=m^{r-1}$
for  all $r$. Consequently, 
$$\lim_{r\to \infty}(b_r^{(p)})^{1/r} = m^{(p)}.$$\bs

The Alexander polynomial of any knot is nonzero, and its
coefficients are relatively prime. Hence the following corollary
is immediate from Theorem 4.2(iii).\bs

\ni {\bf Corollary 4.4.} For any knot $k$ and prime $p$,
$$\lim_{r\to\infty}{(b_r^{(p)}})^{1/r} = 1.$$\bs

Theorem 2.10 and Corollary 4.4 imply that whenever the Alexander
polynomial of
$k$ has Mahler measure greater than $1$, infinitely many
distinct primes occur in the factorization of the torsion numbers
$b_r$. In other words, the homology groups $H_1(M_r,\Z)$ display
nontrivial $p$-torsion for infinitely many primes $p$. Since the
sequence $\{b_r\}$ is a division sequence, the number
of prime factors of
$b_r$ is unbounded.

What about the case in which the Alexander polynomial of $k$ has
Mahler measure equal to $1$? The argument of Section 5.7 of 
[{\bf Go72}] shows that the number of prime factors
remains unbounded as long as the  Alexander polynomial
does not divide
$t^M-1$ for any $M$. If it does divide, then the
torsion numbers
$b_r$ are periodic by Section 5.3 of [{\bf Go72}]
(see also Corollary 2.2 of [{\bf SiWi00}].) Hence we
obtain \bs

\ni {\bf Corollary 4.5.} For any knot, either the torsion numbers
$b_r$ are periodic or else for any $N>0$ there exists an $r$
such that the factorization of $b_r$ has at least $N$ distinct
primes. 
\bs

The proof of Theorem 4.2 requires the following lemma.\bs

\ni {\bf Lemma 4.6.} If
$f(t)=c_0t^n+\cdots+c_{n-1}t+c_n$ is a
nonzero polynomial in $\Z[t]$ with roots $\l_1\cdots,\l_n$ (not
necessarily distinct) in $\Omega_v$, then
$$|c_0|_v\ \prod_{i=1}^n\max\{1, |\l_i|_v \} =|{\rm content}\ f|_v.$$
\bs

\ni {\bf Proof.} The argument that we present is found in
[{\bf LW88}].  Set $a_j= c_j/c_0$ for  $0\le j\le n$, so
$f(t)=c_0(t^n+a_1t^{n-1}+\cdots+a_n)$.  Each $a_j$ is an elementary symmetric function of the
roots
$\l_i$, namely the sum of products of roots taken
$j$ at a time. Using the  ultrametric property
$$|x+y|_v=\max\{|x|_v,|y|_v\},$$
we see that if exactly $k$ values of $|\l_i|_v$ are greater than
$1$, then
$$\max_j |a_j|_v = |a_k|_v = \prod_{j=1}^n\max\{1, |\l_j|_v\}.$$
But
$$\max_j |a_j|_v = \max \{1, |{{c_1}\over
{c_0}}|_v,\ldots, |{{c_n}\over {c_0}}|_v\}={|{\rm content}\ f|_v\over
|c_0|_v}.$$
Hence the lemma is proved. \qed\bs

\ni {\bf Proof of Theorem 4.2.}  In case (i), the pure torsion number
$b_{r_k}$ is equal to 
$\displaystyle |\prod_{\zeta^{r_k}=1}\D(\zeta)|$.
 We have
$$|b_{r_k}|_v = |\prod_{\zeta^{r_k}=1}\D(\zeta)|_v
= |c_0|^{r_k}_v\ \prod_{\zeta^{r_k}=1}\ \prod_{j=1}^n|\zeta-\l_j|_v,$$
where $c_0$ is the leading coefficient of $\D$ and $\l_1, \ldots \l_n$
are its roots. Hence
$$\eqalign {|b_{r_k}|_v^{1/r_k}&= |c_0|_v
\prod_{\zeta^{r_k}=1}\ 
\prod_{j=1}^n|\zeta-\l_j|_v^{1/r_k}\cr
&= |c_0|_v \prod_{j=1}^n \exp \bigr({1 \over
r_k}\sum_{\zeta^{r_k}=1}\log |\zeta -\l_j|\bigr
),}
$$ so that 
$$\lim_{r_k\to \infty}|b_{r_k}|_v^{1/r_k} =
|c_0|_v\prod_{j=1}^n\exp \int_{\T_v}\log|t-\l_j|_v\ d\mu,$$
which by equation (4.1) is equal to 
$$|c_0|_v\ \prod_{j=1}^n \max\{1, |\l_j|_v\}.$$ 
By Lemma 4.6 this is equal to  $|{\rm content}\ \D|_v$. But for 
integers $n$ we have $n^{(p)}= |n|_v^{-1}$.

Now suppose $\cal M$ is cyclic.  As in the proof of Theorem 3.9, we
let $\g$ be the cyclotomic order of $\D$ and consider the subsequence of $b_r$ with
$r$ in a fixed congruence class modulo $\g$. Then starting with Theorem 3.4 we
may apply the argument above with $\D/\Phi$ in place of $\D$ to show
that the limit of $(|b_r|^{(p)})^{1/r}$ along this subsequence is the
$p$-component of the content of $\D/\Phi$.  But content is
multiplicative and cyclotomic polynomials have content 1, so the
limit along all congruence classes is $({\rm content}\ \D)^{(p)}$. 
The result is immediate for direct sums of cyclic modules.

Finally, we can extend the result when $\cal M$ is torsion-free as
an abelian group using Theorem 3.7. But for this case
the content of $\D$ is $1$.
\qed
\bs


\ni {\bf 5. Torsion numbers and links.} A {\it link} is a finite
collection $l=l_1\cup\cdots\cup l_\mu$ of pairwise disjoint knots 
in the $3$-sphere. If a direction is chosen for each component
$l_i$, then the link is {\it oriented}.  Equivalence for links,
possibly oriented, is defined just as for knots. 

The abelianization of the group $G= \pi_1(S^3-l)$ is free abelian of
rank $\mu$ with generators $t_1, \ldots, t_\mu$
corresponding to oriented loops having linking number one with
corresponding components of $l$. When $\mu>1$ there are infinitely
many possible epimorphisms from $G$ to the integers. 

When $l$ is oriented
there is a natural choice for $\chi$, sending each 
generator $t_i$ to $1\in \Z$. In this way we associate to $l$ an
augmented group $(G, \chi)$. As in the special case of a knot, 
$\M$ has a square presentation matrix, and it is
isomorphic to the first homology group of the infinite cyclic
cover of $S^3-l$ corresponding to $\chi$. Again as in the case of a
knot, there is a sequence of $r$-fold cyclic covers $M_r$ of $S^3$
branched over $l$. However, $H_1(M_r;\Z)$ is isomorphic to
$\M/(t^{r-1}+\cdots+t+1)\M$ rather than $\M/(t^r-1)\M$ (see [{\bf
Sa79}]). In the case of a knot the two modules are well
known to be isomorphic (see Remark 5.4(i)).

Motivated by these observations we make the following definitions.
Let $\tM_r$ denote the quotient module $\M/\nu_r\M$, where 
$\nu_r= t^{r-1}+\cdots+t+1$. \bs

\ni {\bf Definition 5.1.} Let $(G, \chi)$ be an augmented group.
The $r$th {\it reduced torsion number} $\tilde b_r$ is the order of
the torsion submodule $T\tM_r$. The $r$th {\it reduced Betti
number} $\tilde \b_r$ is the rank of $\tM$.  \bs

Many results of Section 2 apply to
reduced torsion and Betti numbers with only slight modification. For
example, an argument similar to the proof of Proposition 2.1 shows
that  $\tilde \b_r$ is the number of zeros of the Alexander polynomial
which are roots of unity and different from
$1$, each zero counted as many times as it occurs in the elementary 
divisors $\D_i/\D_{i+1}$; hence $\tilde \b_r$ is periodic
in $r$. Also, when $\tilde \b_r=0$ the reduced torsion number
$\tilde b_r$ is equal to the absolute value of the resultant of 
$\D$ and $\nu_r$. \bs

\ni {\bf Lemma 5.2.} Assume that 
$0\to A\ {\buildrel f \over \to}\ B\ {\buildrel g \over \to}\ C
\to 0$
is an exact sequence of finitely generated abelian groups. If $A$
is finite, then the induced sequence
$$0\to\ A\ {\buildrel f \over \to}\ TB\ {\buildrel g \over \to}\
TC\ \to 0$$
is also exact. \bs

\ni {\bf Proof.} The only thing to check is surjectivity of $g$.
Since the alternating sum of the ranks of $A,B$ and $C$ is zero
and $A$ is finite, the ranks of $B$ and $C$ are equal. By
Lemma 2.3 the homomorphism $g$ maps $TB$ onto
$TC$. \qed\bs

\ni {\bf Proposition 5.3.} Let $(G,\chi)$ be an augmented group such
that $\M$ has a square presentation matrix. If $\D(1)\ne 0$, then for
every $r$, $$\tilde \b_r = \b_r,\quad \tilde b_r =
{{b_r}\over {\d_r}},\eqno(5.1)$$ where $\d_r$ is a divisor of
$|\D(1)|$. Moreover, $\d_{r+\g}=\d_r$, for all $r$, where $\g$ is the
cyclotomic order of $\D$.
\bs

\ni {\bf Proof.} Consider the sequence
$$\M_1\ {\buildrel \nu_r \over \to}\ \M_r\ {\buildrel \pi \over
\to}\ \tilde\M_r\ \to 0,$$
where $\nu_r$ is multiplication by $\nu_r=t^{r-1}+\cdots+t+1$, and
$\pi$ is the natural projection. It is easy to see that the sequence
is exact. From it we obtain the short exact sequence
$$0\to \M_1/{\rm ker}\ \nu_r\ {\buildrel \nu_r \over \to}\ \M_r\
{\buildrel \pi \over \to}\ \tilde\M_r\ \to 0.$$ 
Here $\nu_r$ also denotes the induced quotient homomorphism. 
Since $\D(1)\ne 0$, the module $\M_1$ is finite and hence $\b_r=
\tilde \b_r$. The  order of $\M_1$ is
$|\D(1)|$, and hence the order of $\M_1/{\rm ker}\ \nu_r$ is a divisor
$\d_r$. The second statement of (5.1) follows from Lemmas 5.2 and 3.8.

It remains to show that $\d_r$ has period $\g$. For this let $0\ne a
\in
\M$. The coset $\bar a\in \M_1$ is in the kernel of $\nu_r$ if and 
only if there exists $b\in \M$ such that $\nu_r(a-(t-1)b)=0$.
Clearly this is true if and only if
$\nu_{(\g,r)}(a-(t-1)b)=0$, where $(\g,r)$ denotes the gcd of $\g$
and $r$. Hence the kernel of 
$\nu_r$ is equal to the kernel of $\nu_{(\g,r)}$, and the 
periodicity of $\d_r$ follows.
\qed \bs

\ni {\bf Remarks 5.4.} (i) If $G$ is a knot group, then any two
meridianal generators are conjugate. Consequently $\M_1$ is trivial.
Proposition 5.3 implies that in this case, the torsion numbers $b_r$
and $\tilde b_r$ are equal for every $r$. \ss

\quad (ii) It is well known that for any oriented link $l=l_1\cup l_2$
of two components, $|\D(1)|$ is equal to the absolute value of the
linking number ${\rm Lk}(l_1,l_2)$. (See Theorem 7.3.16 of [{\bf
Ka96}].)
\bs

\ni {\bf Proposition 5.5.} Let $(G,\chi)$ be an augmented group
such that $\cal M$ has a square presentation matrix. Assume that
$\D(t)= (t-1)^q g(t)$, with $g(1)\ne 0$. If $p$ is a prime that does
not divide $g(1)$, then 
$$\tilde \b_{p^k}=0,\quad \tilde b_{p^k}^{(p)}=p^{q k},$$
for every $k\ge 1.$ \bs

The proof of Proposition 5.5 requires: \bs

\ni {\bf Lemma 5.6.} Let $g(t)$ be a polynomial with integer
coefficients, and assume that $p$ is a prime. If $p$ does not divide
$g(1)$, then $p$ does not divide
${\rm Res}(g, t^{p^k}-1)$ for any positive integer $k$. \bs

\ni {\bf Proof of Lemma 5.6.}  Assume that $p$ does not divide
$g(1)$.  Recall that $\Phi_n(t)$ denotes the $n$th cyclotomic
polynomial. From the formula
$$\prod_{d\vert n\atop d>1}\Phi_d(1)=\nu_n(1)=n,$$
we easily derive 
$$\Phi_d(1) = \cases{0&if $d=1$\cr q &if $d=q^k>1$, $q$ 
prime\cr 1& other
$d$.}$$ Consequently, $\Phi_{p^k}$ does not
divide $g$ for any $k>0$, and so  ${\rm Res}(g, t^{p^k}-1)\ne 0$. The
module
${\cal H}=\R1/(g, t^{p^k}-1)$ has order $|{\rm Res}(g, t^{p^k}-1)|$,
and it suffices to prove that ${\cal H}\otimes_\Z \Z/p$ is
trivial. Now, ${\cal H}\otimes_\Z \Z/p$ is
isomorphic to the quotient of the PID $(\Z/p)[t,
t^{-1}]$ by the ideal generated by the greatest
common divisor of $g$ and $t^{p^k}-1$ in this
ring. But $t^{p^k}-1 = (t-1)^{p^k}$ in this
ring,
and $t-1$ does not divide $g$ since $p$ does
not
divide $g(1)$. So the gcd is $1$, and
${\cal H}\otimes_\Z\Z/p$ is trivial.\qed \bs

\ni {\bf Proof of Proposition 5.5.} Let $k$ be any positive
integer. Lemma 5.6 implies that 
${\rm Res}(g, t^{p^k}-1)\ne 0$. Hence $\b_{p^k}$ vanishes, and
therefore
$\tilde \b_{p^k}$ is also zero.
By a result analagous to Proposition 2.5 and the multiplicative
property of resultants
$$\tilde b_{p^k}= |{\rm Res}(\D, \nu_{p^k})| = |{\rm Res}(t-1,
\nu_{p^k})|^q |{\rm Res}(g,
\nu_{p^k})|= (p^k)^q |{\rm Res}(g, \nu_{p^k})|.$$
By Lemma 5.6, $p$ does not divide $|{\rm Res}(g, t^{p^k}-1)|$. Hence
$p$ does not divide
${\rm Res}(g, \nu_{p^k})$, and so
$b_{p^k}^{(p)}=p^{kq}$.\qed\bs

\ni {\bf Corollary 5.7.} (i) Let $M_r$ be the $r$-fold
cyclic cover of $S^3$ branched over a knot. If $r$ is a
prime power $p^k$, then the $p$-torsion submodule of
$H_1(M_r;\Z)$ is trivial. \ss
(ii) Let $M_r$ be the $r$-fold cyclic cover $S^3$ branched over a
$2$-component link $l=l_1\cup l_2$. If $r$ is a power of a prime that
does not divide ${\rm Lk}(l_1,l_2)$, then the $p$-torsion submodule of
$H_1(M_r;\Z)$ is trivial.\bs

\ni {\bf Proof.} Statement (i) was proven in 
[{\bf Go78}]. Here it follows from Proposition 5.5 together
with the well-known fact that $|\D(1)|=1$, whenever $\D$ is the
Alexander polynomial of a knot.
The second statement is a consequence of Proposition 5.5 and Remark
5.4 (ii). \qed\bs

\ni {\bf Proposition 5.8.} Assume that $(G, \chi)$
is an augmented group such that $\M\cong\R1/(\D)$. If
$\D(t)=(t-1)^q g(t)$, where
$g(1)\ne0$, then for every positive integer $r$, there exists 
a positive integer $\d_r'$ such that
$$\tilde b_r = (\d'_r)^q\cdot |T(\R1/(g, \nu_r))|.$$
Moreover, $\d'_{r+\g}=\d'_r$, for all $r$, where $\g$ is the 
cyclotomic order of $\D$. \bs

\ni {\bf Remarks 5.9.} (i) The order $|T(\R1/(g, \nu_r))|$ can 
be found using Proposition 5.3 and Theorem 3.4.

\quad (ii) When $\M$ is a direct sum of cyclic modules,
$\tilde b_r$ can again be found by applying Proposition 5.5 to each
summand. When $\M$ is not a direct sum of cyclic modules but is
torsion free as an abelian group, a result analagous to Theorem 3.7
can be found by replacing $t^r-1$ everywhere by $\nu_r$ in the proof.
As in Section 3, the torsion numbers $\tilde b_r$ are then seen to
satisfy a linear homogeneous recurrence relation. 
\bs

\ni {\bf Proof of Proposition 5.8.} Consider the exact sequence
$$0\to {\rm ker}\ g \to \R1/((t-1)^q,\nu_r)\ {\buildrel g \over \to}\
\R1/((t-1)^qg,
\nu_r)\ {\buildrel \pi \over \to}\ \R1/(g, \nu_r)\to 0, $$
where the first homomorphism is inclusion, the second 
is multiplication by
$g$, and the third is the natural projection. The order of 
$\R1/((t-1)^q,\nu_r)$ is equal to $|{\rm Res}((t-1)^q, \nu_r)|$,
which is equal to $r^q$. The kernel of $g$ is generated by
$\nu_r/f_r$, where
$f_r$ is the greatest common divisor of
$g$ and $\nu_r$. Notice that $f_{r+\g}=f_r$, for all $r$. Lemmas
5.2 and 3.8 complete the proof. \qed\bs

We conclude with a generalization of Corollary 5.7(ii). 

When $(G,
\chi)$ is the augmented group corresponding to a $2$-component link
$l$, the epimorphism
$\chi$ factors through $\eta: G \to G_{ab} \cong
\Z^2$.  For any finite-index subgroup $\L \subset
\Z^2$ there is a
$|\Z^2/\L|$-fold cover of $S^3$ branched over $l$
corresponding to the map $G\to \Z^2\to \Z^2/\L$.
The cover
$M_r$ is a special case corresponding to the
subgroup
$\L$ generated by
$t_1-t_2, t_1^r, t_2^r$. We denote the rank of $H_1(M_\L;\Z)$ by
$\b_\L$ and the order  $|TH_1(M_\L;\Z)|$ by $b_\L$.
\bs

\ni {\bf Theorem 5.11.} Let $l=l_1\cup l_2$ be a link in $S^3$. If
$p$ is a prime that does not divide ${\rm Lk} (l_1,l_2)$, then
$\b_\L=0$ and
$b_\L^{(p)}=1$ for any subgroup $\L\subset \Z^2$ of index $p^k,
k\ge 1$. \bs

\ni {\bf Proof.} Let $\M_\eta$ be the kernel of $\eta$. We consider the
dual $\M_\eta^\wedge$, which is a compact abelian group with a
$\Z^2$-action by automorphisms induced by conjugation in $G$ by $t_1$
and $t_2$. The automorphism induced by ${\bf n} \in \Z^2$ is denoted
by $\s_{\bf n}$; the automorphims induced by $(1,0)$ and $(0,1)$ are
abbreviated by $\s_1$ and $\s_2$, respectively. The dual
$\M_\eta^\wedge$ can be identified with a subspace of ${\rm
Fix}_\L(\s)= \{\rho\in \M_\eta^\wedge\ :\ \s_{\bf n}\rho=\rho\ {\rm
for\ all}\ {\bf n}\in \L\}$. Details can be found in [{\bf SW00}].

From the elementary ideals of $\M_\eta$ a sequence of $2$-variable
Alexander polynomials $\D_i(t_1,t_2)$ is defined; when $i=0$, setting
$t_1=t_2=t$ recovers $\D(t)$. By [{\bf Cr65}],
$\D_0(t_1, t_2)$ annihilates
$\M_\eta$. Hence $\D_0(\s_1, \s_2)\rho =0$ for all $\rho\in
\M_\eta^\wedge$. Consequently, if $\s_{\bf
n}\rho=\rho$ for all ${\bf n}\in \Z^2$ then 
$0=\D_0(\s_1,\s_2)\rho= \D_0(1,1)\rho=
\D(1)\rho.$ Recall that $\D(1)={\rm Lk}(l_1,l_2)$. 

Let
$$Y= \{\rho: \M_\eta \to \Z/p\ :\ \s_{\bf n}\rho=\rho\ {\rm for\ all}\
{\bf n}\in \L\}.$$
We identify $\Z/p$ with the group of $p$th roots
of unity, so that $Y$ is contained in
$\M_\eta^\wedge$. It is a subspace of
${\rm Fix}_\L(\s)$ invariant under the
$\Z^2$-action, and it contains a subspace isomorphic to
$\M_\eta\otimes_\Z\Z/p$. It suffices to prove that $Y$ is trivial. 

Our hypothesis that $p$ does not divide
the linking number of $l_1$ and $l_2$ implies that $\D_0(t_1, t_2)$ is
not zero. Consequently, $Y$ is a finite $p$-group and so its order is
a power of $p$. In view of the second paragraph,
the hypothesis also implies that the only point
fixed by the
$\Z^2$-action is trivial. But 
$$|Y| = \sum |{\cal  O}_\rho| = \sum |\Z^d/{\rm stab}(\rho)|,$$
where the sums are taken over distinct orbits
${\cal O}_\rho$ and stabilizers ${\rm stab}(\rho)$, respectively.
Each stabilizer contains $\L$, and so $|\Z^d/{\rm stab}(\rho)|$ is
a divisor of $p^k$ whenever $\rho \ne 0$. Hence
$|Y|$ is congruent to $1$ mod $p$. Since $|Y|$
is a power of $p$, the
subspace
$Y$ must be trivial. \qed\bs

\centerline{\bf References.} 
\bs
\baselineskip=12 pt

\item{[{\bf Al28}]} J.~W. Alexander, ``Topological invariants of
knots and links,'' {\sl Trans.\ Amer.\ Math.\ Soc.\ \bf 30} (1928),
275--306. 
\ss
\item{[{\bf AB27}]} J.~W. Alexander and G.B.~Briggs, ``On types of
knotted curves,'' {\sl Ann.\ of\ Math.\ \bf 28} (1927), 562--586.
\ss
\item{[{\bf Ba77}]} A. Baker, ``The theory of linear forms
in logarithms,'' in Transcendence Theory: Advances and 
Applications (Proc., Univ. Cambridge, Cambridge,
1976) Academic Press, London, 1977.
\ss
\item{[{\bf Bo81}]} D.~W. Boyd, ``Speculations concerning the
range of Mahler's measure,'' {\sl Canad.\ Math.\ Bull.\ \bf
24} (1981), 453--469.
\ss

\item{[{\bf Br92}]} R.~A. Brualdi, Introductory
Combinatorics, Second Edition, Prentice Hall, NJ 1992.
\ss
\item{[{\bf Cr65}]} R.~H. Crowell, ``Torsion in link modules,''
{\sl J.\ Math.\ Mech.\ \bf 14} (1965), 289--298.
\ss
\item{[{\bf EEW00}]} M. Einsiedler, G.~R. Everest and T. Ward,
``Primes in sequences associated to polynomials (after Lehmer),''
{\sl London\ Math.\ Soc.\ Journ.\ Comput.\ Math.\ \bf 3} (2000),
125--139.
\ss
\item{[{\bf Ev99}]} G.R.~Everest, ``On the elliptic analogue of
Jensen's formula,'' {\sl J.\ London\ Math.\ Soc.\ \bf 59} (1999),
21--36.
\ss

\item{[{\bf EF96}]} G.R.~Everest and B.N\'i Fhlath\'uin, ``The
elliptic Mahler measure,'' {\sl Math.\ Proc.\ Camb.\ Phil.\ Soc.\ 
\bf 120} (1996), 13--25.
\ss
\item{[{\bf EW99}]} G. Everest and T. Ward, Heights of polynomials
and entropy in algebraic dynamics, Springer-Verlag, London, 1999. 
\ss
\item{[{\bf Fo54}]} R.~H. Fox, ``Free differential calculus. III
Subgroups,'' {\sl Annals\ of\ Math.\ \bf 59} (1954), 196--210.
\ss
\item{[{\bf Ge35}]} A.~O. Gelfond, ``On the
approximation of transcendental numbers by
algebraic numbers,'' {\sl Dokl.\ Akad.\ Navk.\
SSSR\ \bf 2} (1935), 177--182.
\ss
\item{[{\bf Go72}]} C.~McA. Gordon, ``Knots whose branched coverings
have periodic homology,'' {\sl Trans.\ Amer.\ Math.\ Soc.\ \bf 168}
(1972), 357--370.
\ss
\item{[{\bf Go78}]} C.~McA. Gordon, ''Some aspects of classical knot
theory,'' in Knot Theory (Proc., Plans-sur-Bex, Switzerland,
1977), Lecture Notes in Mathematics {\bf 685} (J.~C. Hausmann),
Springer-Verlag, Berlin, Heidleberg, New York, 1978, 1--60.
\ss
\item{[{\bf GS91}]} F. Gonz\~alez-Acu\~na and H. Short, ``Cyclic
branched coverings of knots and  homology spheres,'' {\sl Revista\
Math.\ \bf 4} (1991), 97--120.
\ss
\item{[{\bf Hi81}]} J.~A. Hillman, Alexander ideals of l

inks,
Lecture Notes in Math. {\bf 895}, Springer-Verlag, Berlin,
Heidelberg, New York, 1981.
\ss
\item{[{\bf HS97}]} J.~A. Hillman and M. Sakuma, ``On the homology of
links,'' {\sl Canad.\ Math.\ Bull.\ \bf97} (1997), 309--315.
\ss
\item{[{\bf Ka96}]} A. Kawauchi, A Survey of Knot Theory,
Birkh\"auser, Basel, 1996. 
\ss
\item{[{\bf La65}]} S. Lang, Algebra, Addison-Wesley, Reading, 1971.
\ss
\item{[{\bf Le33}]} D.~H. Lehmer, ``Factorization of certain
cyclotomic functions,'' {\sl Annals\ of\ Math.\ \bf 34} (1933), 461
-- 479. 
\ss
\item{[{\bf Li97}]} W.~B. Lickorish, An Introduction to Knot Theory,
Springer-Verlag, Berlin, 1997. 
\ss
\item{[{\bf LW88}]} D. Lind and T. Ward, ``Automorphisms of
solenoids and $p$-adic entropy,'' {\sl Ergod.\ Th.\ and\ Dyn.\
Syst.\ \bf 8} (1988), 411 -- 419.
\ss
\item{[{\bf LSW90}]} D. Lind, K. Schmidt and T. Ward, ``Mahler
measure and entropy for commuting automorphisms of compact
groups,'' {\sl Invent.\ Math.\ \bf 101} (1990), 593--629.
\ss
\item{[{\bf MM82}]} J.~P. Mayberry and K. Murasugi,
``Torsion-groups of abelian coverings of links,'' {\sl Trans.\
Amer.\ Math.\ Soc.\ \bf 271} (1982), 143--173.
\ss
\item{[{\bf Me80}]} M.~L. Mehta, ``On a relation between
torsion numbers and Alexander matrix of a knot,'' {\sl Bull.\
Soc.\ Math.\ France\ \bf 108} (1980), 81--94.
\ss
\item{[{\bf Ne65}]} L.~P. Neuwirth, Knot Groups, Princeton Univ.\
Press, Princeton, NJ, 1965.
\ss
\item{[{\bf Ri90}]} R. Riley, ``Growth of order of homology of
cyclic branched covers of knots,'' {\sl Bull.\ London\ Math.\ Soc.\
\bf 22} (1990), 287--297.
\ss
\item{[{\bf Ro76}]} D. Rolfsen, Knots and Links, Publish or Perish, 
Berkeley, CA, 1976. 
\ss
\item{[{\bf Sa79}]} M. Sakuma, ``The homology groups of abelian
coverings of links,'' {\sl Math.\ Sem.\ Notes\ Kobe\ Univ.\ \bf 7}
(1979), 515--530.
\ss
\item{[{\bf Sa95}]} M. Sakuma, ``Homology of abelian coverings
of links and spatial graphs,'' {\sl Canad.\ J.\ Math.\ \bf 47}
(1995), 201--224. 
\ss
\item{[{\bf Sc95}]} K. Schmidt, Dynamical Systems of Algebraic 
Origin, Birkh\"auser Verlag, Basel, 1995.
\ss
\item{[{\bf SW00}]} D.~S. Silver and S.~G. Williams, ``Mahler
measure, links and homology growth,'' {\sl Topology}, in press. Preprint on
Mathematics ArXive, GT/00003127.
\ss
\item{[{\bf St00}]} W.~H. Stevens, ``Recursion formulas for some
abelian knot invariants,'' {\sl Journal\ of\ Knot\ Theory\ and\ its\
Ramifications\ \bf 9} (2000),  413--422.
\ss
\item{[{\bf We80}]} C. Weber, ``Sur une formule de R.~H. Fox
concernant l'homologie d'une rev\^etment ramifi\'e,'' {\sl
Enseignement\ Math.\ \bf 25} (1980), 261--272. 
\ss
\item{[{\bf Yo86}]} R.~M. Young, ``On 	Jensen's formula and 
$\int_0^{2 \pi}\log |1-e^{i\theta}|d\theta$,'' {\sl Amer.\ Math.\
Monthly\ \bf93} (1986), 44--45.
\ss
\item{[{\bf Za32}]} O. Zariski, ``On the topology of
algebroid singularities,'' Amer. J. Math. {\bf 54}
(1932), 453--465. 
\bs
\item{} Dept.\ of Mathematics and Statistics, Univ.\ of South Alabama, Mobile, AL 
36688-0002 e-mail: silver@mathstat.usouthal.edu, williams@mathstat.usouthal.edu

\end
\ni {\bf 4. On a theorem of Fried.} Let $\D$ be a  polynomial with
real coefficients that is {\it reciprocal} in the sense that $t^{{\rm
deg}\ \D}\cdot \D(t^{-1})=\D(t)$.  In [{\bf Fr88}] D. Fried proved
that the sequence of resultants
${\rm Res} (\D, t^r-1)$ determines $\D$ up to multiplication by a
unit in $\R1$, provided that
$\D$ is reciprocal and no zero of
$\D$ is a root of unity. Much of
the motivation for Fried's result comes from situations in which $\D$
has integer coefficients.

For any polynomial $\D$ with integer coefficients, we can consider the
module $\M=\R1/(\D)$. It arises from an augmented group that we need
not  explicitly describe. For any $r\ge 1$, the $r$th Betti number
$\b_r$ and the $r$th torsion number $b_r$ can be computed directly
from $\D$. As in Section 2, $\b_r$ is the number of distinct zeros of
$\D$ that are $r$th roots of unity, while $b_r$ is given by formula
(3.2).

The sequences $b_r$ and $\b_r$ determine ${\rm Res} (\D,
t^r-1)$: If $\b_r\ne0$, then ${\rm Res} (\D,
t^r-1)=b_r$; if $\b_r=0$ then ${\rm Res} (\D, t^r-1)=0$. We
prove that $b_r$ and $\b_r$ together determine $\D$, given some conditions on
$\D$.
\bs

\ni {\bf Theorem 4.1.} Assume that $\D$ is a reciprocal polynomial
with integer coefficients. Suppose $\D$ is not divisible by $t-1$ and has no root
on the unit circle that is not a root of unity. Then the sequences
$\{b_r\}$ and
$\{\b_r\}$ of torsion and Betti numbers, respectively, 
associated to $\M=\R1/(\D)$ together determine $\D$ up to
multiplication by a unit.
\bs

\ni {\bf Remarks 4.2.} 1. Betti numbers are fairly crude invariants,
and they alone do not suffice to determine $\D$. For example, the
reciprocal polynomials $(t-1)^2$ and $(t-1)^4$ share the sequence
$1,1,1,\ldots$ of Betti numbers but have 
torsion number sequences $\{r\}$, $\{r^3\}$ respectively.

Torsion numbers are more sensitive. Yet, as with Betti numbers, they
alone do not suffice to determine $\D$. The polynomials
$t-1$ and $t^2-1$, for example, each have a trivial sequence
$1,1,\ldots$ of torsion numbers but different sequences of Betti
numbers. Reciprocal examples are given in [{\bf Fr88}].

2. As in [{\bf Fr88}] the reciprocality of
$\D$ is needed for the theorem to hold. For example, the polynomials
$f(t)=t^2+t+2$ and $g(t)= 2t^2+t+1$ have the same sequences of torsion
and Betti numbers. Notice that $f(t^{-1})$ and $g(t)$ are equal up
to a unit. Possibly this is the only indeterminancy that remains if
the hypothesis of reciprocity is dropped. 

3. We do not know it we can drop the hypotheses that $\D$ is not divisible by 
$t-1$ and has no root
on the unit circle that is not a root of unity. We will see in the proof that we
can detect from the sequences $b_r$ and $\b_r$ whether these conditions are
satisfied.

4. We recall from the proof of Theorem 2.10 that $b_r$ is
the number of connected components of ${\rm Fix}\ \s^r$ acting on
the dual space $\M^{\wedge}$. Each connected component is a torus 
of dimension $\b_r$. Theorem 4.1 implies that the numbers of connected
components together with their dimensions determine $\D$ up to
multiplication by a unit, and hence they determine all of the dynamics
of $\s$.
\bs

\ni {\bf Proof of Theorem 4.1.} Write $\D$ as
$\Phi_{n_1}^{e(n_1)}\cdots \Phi_{n_k}^{e(n_k)}g$, where 
$\Phi_{n_1}, \dots,\Phi_{n_k}$  are distinct cyclotomic
polynomials and $g$ is a reciprocal polynomial with no root on the unit circle
As a first step we prove that the sequence of 
Betti numbers $\b_r$ determines $n_1, \ldots, n_k$.
For any positive integer $d$, define 
$$\hat \phi(d)= \cases{\phi(d) & if $e(d)\ge 1$, \cr
0& if $e(d)=0$}$$
(as before, $\phi$ is Euler's phi function).
Then 
$$\b_r= \sum_{d\vert r}\hat\phi(d).$$
By the M\"obius inversion formula,
$$\hat\phi(r)=\sum_{d\vert r}\mu(d)\b_{r/d}.$$
Here $\mu$ is the M\"obius $\mu$-function, defined by 
$$\mu(d)=\cases{1&if $d=1$,\cr
                 0&if $p^2\vert d$ for some prime $p$,\cr
                 (-1)^q&if $d$ is the product of $q$ distinct
                 primes.}$$
Using the sequence $\{\b_r\}$ we determine $\{\hat\phi(r)\}$. Then
$\Phi_n$ divides $\D$ if an only if $\hat\phi(n)>0$. Thus we find
$n_1,\ldots, n_k$. Note that $\b_1=0$ if and only all $n_j>1$.

As a second step we use $\{b_r\}$ to find
the exponents $e(n_1),\ldots, e(n_k)$. Select a prime $p$ that does
not divide any $n_j$ or $b_{n_j}$, for $1\le j\le k$. For any integer $n$ we denote
by $n^{(p)}$ the highest power of $p$ dividing $n$. Consider a particular $n_i$. Let
$\Phi=\prod_{n_j\vert n_i}\Phi_{n_j}$. By Theorem 3.4 and the multiplicativity of resultants,
$b_{n_ip}^{(p)}$ is equal to
$$\prod_j |{\rm Res}(\Phi_{n_j}^{e'(n_j)},
(t^{n_ip}-1)/\Phi)|^{(p)}|{\rm Res}(g,(t^{n_ip}-1)/\Phi)|^{(p)},$$
where 
$$e'(n_j)=\cases{e(n_j)-1&if $n_j\vert\ n_i$,\cr
                 e(n_j)&if $n_j\not\vert n_i$.}$$
The following fact about resultants of cyclotomic
polynomials is due to T. Apostol (see Theorems 1, 3 and 4 of [{\bf
Ap70}]). For
$1\le n<m$,
$${\rm Res}(\Phi_n, \Phi_m)=\cases{p^{\phi(n)}&if $m=np^i$, for
some prime power $p^i$,\cr 1& otherwise,}$$ where $\phi$ denotes
Euler's phi function.  Using this together with the factorization
$$t^{n_ip}-1 = \prod_{m\vert n_ip}\Phi_m,$$
we have 
$$b_{n_ip}^{(p)} = \prod_{n_j\vert n_i} p^{e'(n_j)\phi(n_j)}
\biggr|{\rm Res}(g,{t^{n_ip}-1 \over \Phi})\biggr|^{(p)}.$$

We claim that $|{\rm Res}(g, (t^{n_ip}-1)/\Phi)|^{(p)} = 1,$ for
almost all primes $p$. Since no zero of $g$ is a root of unity, 
$|{\rm Res}(g,(t^{n_ip}-1)/\Phi)|=|\R1/(g, (t^{n_ip}-1)/\Phi)|$. In
order to  prove the claim it suffices to show that the tensored module
$\R1/(g, (t^{n_ip}-1)/\Phi)\otimes_\Z\Z/p$ is trivial. For this, note
that $t^{n_ip}-1\equiv (t^{n_i}-1)^p$ mod $p$. Hence 
$\R1/(g,(t^{n_ip}-1)/\Phi)\otimes_\Z\Z/p \cong 
\R1/(g, (t^{n_i}-1)^p/\Phi)\otimes_\Z\Z/p$. But 
$$|\R1/(g, (t^{n_i}-1)^p/\Phi)| = |\R1/(g, (t^{n_i}-1/\Phi)^p||\R1/(g,
\Phi^{p-1})|$$ $$= |{\rm Res}(g, (t^{n_1}-1)/\Phi)|^p|{\rm
Res}(g,\Phi)|^{p-1},$$ 
which divides $b_{n_i}^p\cdot|{\rm Res}(g,\Phi)|^{p-1}.$ We have
already ensured by our choice of $p$ that $p$ does not divide
$b_{n_i}$. Only finitely many primes divide ${\rm Res}(g, \Phi)$. 
Hence the above claim is established. 

We choose a prime $p$ such that $\log_p b_{n_ip}^{(p)}$ is a minimum,
for each $i=1,\ldots, k$. From what has been said, this minimum value
is $\sum_{n_j\vert n_i}(e(n_j)-1)\phi(n_j)$. 
We easily recover the exponents $\{e(n_j)\}$, beginning with the $n_j$
that are minimal under the partial ordering of divisibility.  (This is
M\"obius inversion on the finite partially ordered set $\{n_1,\ldots,n_k\}$.)

The third step of the proof is to find the sequence $\{|{\rm Res}(g,
t^r-1)|\}$. Define 
$$b_r^\circ = \prod_{\zeta^r=1\atop \D(\zeta)\ne 0}|g(\zeta)|,\quad
d_r=\prod_{\Phi(\zeta)=0}|g(\zeta)|,$$
where as before $\Phi= \prod\{\Phi_d : d\vert r, \Phi_d\vert \D\}$.
Clearly, $|{\rm Res}(g, t^r-1)|=b_r^\circ d_r$ for each $r$. Our 
choice of prime $p$ ensures that
$d_{rp^i}=d_r$, for all $i\ge1$.  In view of steps 1 and 2 and Theorem 3.4, we can
find $\{b_r^\circ\}$. We compute
$\{d_r\}$ as follows. 

Write $g(t)=c_0\prod_j(t-\a_j)$, where the $\a_j$ are the zeros of
$g$ (not necessarily distinct). Then 
$$b_r^\circ d_r = |{\rm Res}(t^r-1,g)| =
|c_0|^r\prod_j|1-\a_j^r|.$$   
Consider those $n_1, \ldots, n_k$ that
do not divide $r$, and let $q$ be the product of the distinct primes
that divide at least one of them but do not divide $r$. If there is
no such prime, then let $q=1$. It is easy to see that no $n_j$, $1\leq j\leq k$,
that is greater than one divides $rp^i+q$ for any $i\geq 1$. (Consider primes
divisors of
$n_1,
\ldots, n_k$.) Hence $d_{rp^i+q}=d_1$, for any $i\ge 1$, and so 
$${b_{rp^i+q}^\circ\over {b_{rp^i}^\circ}} = {d_r\over d_1}|c_0|^q \prod_j
\biggr\vert {{1-\a_j^{rp^i+q}}\over {1-\a_j^{rp^i}}}\biggr\vert.$$
As $i$ goes to infinity, this ratio approaches 
$${d_r\over d_1}|c_0|^q \prod_j {\rm max}(1, |\a_j|^q) = {d_r\over d_1} M^q,$$
where $M$ is the Mahler measure of $\D$. (If $g$ has a root on the unit circle that
is not a root of unity, this sequence will diverge [REF]). We find
$M$ by taking
$r=1$. If no
$n_j=1$ then $d_1=1$, so we can find the remaining $d_r$. 

Having found $\{|{\rm Res}(g, t^r-1)|\}$, we apply the main result
of [{\bf Fr88}] to recover the polynomial $g$. Since $n_1, \ldots,
n_k$ and $e(n_1), \ldots, e(n_k)$ have been found, the polynomial
$\D$ is determined. \qed\bs

\bs

\bs